\documentclass{article}
\usepackage{CJK,CJKnumb,CJKulem,times,dsfont,ifthen,mathrsfs,latexsym,amsfonts,color}
\usepackage{amsmath,amsthm,makeidx,fontenc,amssymb,bm,graphicx,psfrag,listings,curves,extarrows,mathrsfs}
\usepackage{xypic}
\usepackage[all]{xy}
\newtheorem{theo}{Theorem}[section]
\newtheorem{pro}[theo]{Proposition}

\newtheorem{con}[theo]{Conjecture}
\theoremstyle{definition}

\numberwithin{equation}{section}

\begin{document}
\begin{CJK*}{GBK}{song}
\title{\bf   A minimax  argument to a stronger version of the  Jacobian conjecture\thanks{Supported by NSFC.\quad\quad
E-mails: liuw16@mails.tsinghua.edu.cn }}

\date{}
\author{
{\bf Wei Liu}      \\
\footnotesize  \it  Department of Mathematical Sciences, Tsinghua University,\\
\footnotesize \it  Beijing 100084, China \\\\}
 

\maketitle

\vskip0.6in

\begin{center}
\begin{minipage}{120mm}
\begin{center}{\bf Abstract}\end{center}
The main result of this paper is to prove the strong real Jacobian conjecture under the symmetric assumption (Theorem 1.6) and reveals the link between it and the Jacobian conjecture (Proposition 1.3).  Precisely,
we assume that  $F: \mathbb{R}^n \to \mathbb{R}^n$  is of  $C^1$ map, $n\geqslant 2$, if for some $\varepsilon >0$,
	$ 0\notin Spec(F)~~\mbox{and}~ Spec(F+F^T) \subseteq  (-\infty,-\varepsilon)~\mbox{or} ~(\varepsilon,+\infty),$
	where $Spec (F)$ denotes all  eigenvalues of $JF$ and $Spec (F+F^T)$ denotes all eigenvalues of $JF+JF^T$, then we show that  $F$ is  injective.
 It is proved by  using a  minimax  argument.

 \vskip0.23in

{\it   Key  words:}  Jacobian Conjecture; Minimax method;   Injective.

 \vskip0.23in

 MSC(2010): 14R15; 35A15; 18G05.
\vskip0.23in

\end{minipage}
\end{center}
\vskip0.26in
\newpage
\section{Introduction}

In 1939, Keller (see \cite{Kel}) stated the following Conjecture:
\begin{con}\label{C1.1}{\bf(Jacobian Conjecture)}.
	Let $  F: \mathit{k}^n \to \mathit{k}^n$ be a polynomial map, where $ \mathit{k} $  is a field of characteristic 0. If the determinant for its jacobian of the polynomial map is a non-zero constant, i.e., $ \det JF(x) \equiv C \in \mathit{k}^{*},~   \forall x\in \mathit{k}^n$,  then $ F(x) $ has a polynomial inverse map.
\end{con}

This is a long-standing conjecture, it is still open even in the case $n=2$. There are many partial results on this conjecture, see for example \cite{BCW} and \cite{Ess1}.  Chamberland and  Meisters discovered a sufficient condition for Conjecture \ref{C1.1} and formulated the following  conjecture:
\begin{con}(\cite{CM},~Conjecture 2.1)\label{C1.4}
Let $F:\mathbb{R}^n\rightarrow \mathbb{R}^n$ be a ${C^1 }$ map.
Suppose there exists  an $\varepsilon>0$ such that $|\lambda|\geqslant\varepsilon$~ for all the eigenvalues $\lambda$ of $F'(x)$ and all $x\in \mathbb{R}^n$.
Then $F$ is injective.
\end{con}

In this note, we will show that the following conclusion is true:

\begin{pro}\label{pro}
Conjecture \ref{C1.4} implies Conjecture \ref{C1.1}. That is,
if Conjecture \ref{C1.4} is true for all $n$, then  Conjecture \ref{C1.1} is also true.
\end{pro}

In the current paper, let $Spec(F)$ denote  the set of all (complex) eigenvalues  of  $JF(x)$ and let  $Spec(F+F^T)$  denote
the set of all (complex) eigenvalues of $JF+\{JF\}^T$, for any $x\in \mathbb{R}^n$,
where $JF$ is  the Jacobian matrix of $F$.
 Fernandes,  Gutierrez, and  Rabanal proved Conjecture \ref{C1.4}
in dimension $n=2$ and obtained the following strong  theorem.

\begin{theo}\label{T1}(\cite{FGR})
	Let $F=(f,g):\mathbb{R}^2\rightarrow \mathbb{R}^2$ be a differentiable map.
	For some $\varepsilon >0$,
	if
	\begin{equation}\label{S1}
	Spec (F)\cap [0,\varepsilon )=\varnothing,
	\end{equation}
	then $F$ is injective.
\end{theo}

Recently, Liu and Xu directly prove the Conjecture \ref{C1.4} under an addtional condition and have the following theorem.

\begin{theo} \label{t1.4}(\cite{LX})
Suppose that $F: \mathbb{R}^n\rightarrow \mathbb{R}^n$ is a $C^1$ map, $n\geqslant 2$.
If there exists $\varepsilon>0$, such that
 \[ Spec(F)\in \mathbb{C} \backslash (-\varepsilon,\varepsilon)\times (-i\varepsilon,i\varepsilon)~~~\mbox{and}~~~ Spec(F+F^T) \subseteq  (-\infty,-\varepsilon)~\mbox{or} ~(\varepsilon,+\infty),\]
then $F$ is injective.
\end{theo}

In the current paper, we can obtain the following stronger results.
\begin{theo} \label{main}
	Let $F: \mathbb{R}^n \to \mathbb{R}^n$  be a $C^1$ map,  $n\geqslant 2$.
	For some $\varepsilon >0$, if
	\[ 0\notin Spec(F)~~~\mbox{and}~~~ Spec(F+F^T) \subseteq  (-\infty,-\varepsilon)~\mbox{or} ~(\varepsilon,+\infty),\]
	then $F$ is  injective.
\end{theo}

Considering the Pinchuk's counterexample (see \cite{Pin}) that the eigenvalues sometimes can tend to zero
although they are positive, we have the following result under adding  the conditon of
Spec$(F+F^T)$,

\begin{theo} \label{main2}
Let $F$ be a $C^1$ map from $\mathbb{R}^n$ to $\mathbb{R}^n$.
Suppose that $0\notin Spec(F)$
and that 
one of the following conditions holds:

(i) $Spec(F+F^T)\subseteq (-\infty,0)$ and $\exists $
$M _1,~ M _2> 0$ such that
\[\sum_{\lambda\in Spec(F+F^T) }\lambda>-M_1 ~\mbox{and}~
\prod_{\lambda\in Spec(F+F^T)}|\lambda |>M_2;\]

(ii) $Spec(F+F^T)\subseteq (0,+\infty)$and $\exists $
$M_1,~ M_2> 0$ such that
\[\sum_{\lambda\in Spec(F+F^T) }\lambda<M_1 ~\mbox{and}~
\prod_{\lambda\in Spec(F+F^T)}\lambda >M_2.\]

Then $F$ is injective.
\end{theo}

\vskip0.09in

\section{Minimax Method}

We firstly recall the following Mountain Pass theorem.
\begin{theo}\label{mM} Let $X$ be a Banach space, and  $I\in C^1(X, \mathbb{R})$. Let $ \Omega  \subset X$ be an open set with $u_0 \in \Omega $
	and $u_1 \notin  \Omega .$  Set
	\[\Gamma  = \left\{ {\left. {\gamma \in C([0,1],X)} \right|\gamma (i) = {u_i},i = 0,1} \right\}\]
	and
	\begin{equation}
	c = \mathop {\inf }\limits_{\gamma  \in \Gamma } \mathop {\sup }\limits_{t \in [0,1]} I\big(\gamma (t)\big).
	\end{equation}
	If further 
\begin{itemize}
	\item[(a)] $\alpha  = \mathop {\inf }\limits_{\partial \Omega } I(u) > \max \{ I({u_0}),I({u_1})\}$;
	\item[(b)] $I$ satisfies $(PS)_c$ condition.
\end{itemize}
Then $c$ is a critical value of $I$.
\end{theo}


\noindent{\bf Proof of Proposition \ref{pro}.} We show that Conjecture \ref{C1.4} implies  Conjecture \ref{C1.1}. Firstly, it is sufficient to consider the case $k=\mathbb{C}$ for Conjecture \ref{C1.1} by Lefschetz principle (see \cite{Ess}). Next, it  further reduces to prove that for $k=\mathbb{C}$, $F$ is injective  in conclusion of Conjecture \ref{C1.1} (see \cite{CR}). Finally, we prove $F$ is injective for Conjecture \ref{C1.1} in $\mathbb{C}$ by Conjecture \ref{C1.4}.

Let $F$ be a polynomial map from ${\mathbb{C}^{n}}$ to ${\mathbb{C}^{n}}$ defined as following:
$$({x_1},x_2,...,{x_n}) \to ({F_1},F_2,...,{F_n}).$$
Define a polynomial map $\overline F:{\mathbb{R}^{2n}} \to {\mathbb{R}^{2n}} $ as 
\[\left( {{\mathop{\rm Re}\nolimits} {x_1},{\mathop{\rm Im}\nolimits} {x_1},...,{\mathop{\rm Re}\nolimits} {x_n},{\mathop{\rm Im}\nolimits} {x_n}} \right) \to \left( {{\mathop{\rm Re}\nolimits} {F_1},{\mathop{\rm Im}\nolimits} {F_1},...,{\mathop{\rm Re}\nolimits} {F_n},{\mathop{\rm Im}\nolimits} {F_n}} \right).\]
Then $\det J\overline F=|\det JF|^2$. Consequently, $\det JF$ is not zero complex constant if and only if $\det J \overline F$ is
not zero real constant.
\vskip0.12in
For $\overline F$,  we consider  $\overline F(x)=x-H(x)$, where $H$ is a cube-homogeneous polynomial map and
$JH(x)$ is a nilpotent matrix, i.e., $JH(x)^{2n}=0~$ (see \cite{BCW}). From linear algebra, we know that $JH(x)$ is nilpotent if and only if its characteristic polynomial $\det \big(\mu I - JH(x)\big)=\mu^{2n}$. Therefore, we let $\mu=\lambda -1$ and compute the characteristic polynomial of $J \overline F(x)$,
\[\det \big(\lambda I-J\overline F(x)\big) = \det\big( (\lambda-1)I - JH(x) \big) =\det\big( \mu I - JH(x) \big)=\mu^{2n}=0.\]
\vskip0.12in
Thus, $\mu=\lambda -1=0$. That is, $Spec(\overline F)=\{1\}$. It implies that $Spec(\overline F)$ can not tend to zero.
By Conjecture \ref{C1.4}, we obtain that 
$\overline F$ is injective. Oboviously, the map $F$ is injective if and only if $\overline F$ is injective. Thus, $F$
is injective. By the known result (see \cite{CR}), we have $F$ is bijective. Furthermore the
inverse is also a polynomial map.
 \hfill $\Box$


\section{The proof of Theorem \ref{main}}
\begin{proof}
 Suppose by contradiction that $F$ is not injective, then $F(a)=F(b)$ for some
$a, b\in \mathbb{R}^n, a \ne  b$.
We define $G\left( X \right) = F\left( {X + a} \right) - F\left( {b} \right), \forall X \in \mathbb{R}^n$. Then $G(0)=0$ and putting $c=b-a$, we have $c\ne 0$ and $G(c)=0$.
Let
$I\left( X \right) = G\left( X\right)^T G{\left( X \right)}, ~\forall X \in \mathbb{R}^n.$
Thus $I'(X)=2G(X)^TG'(X)$ and $I(c)=I(0)=0$.

Observe $G'(X)=F'(X+a)$, so $G'(X)$ has no zero eigenvalue.
Therefore,
\[\det G'(X) \ne 0,~\forall X \in \mathbb{R}^n.\]
If $I'(X) = 0,~\forall X \in \mathbb{R}^n $,
i.e., $G(X)^TG'(X)=0,~\forall X \in \mathbb{R}^n $, then $G'(X)G(X)=0,~\forall X \in \mathbb{R}^n.$
So $G(X)=0$ and $I(X)=0$.  Next, we prove that $I(X)$ satisfies the geometric condition-(a) in Theorem 2.3.
Since $I(c)=I(0)=0$,  it is sufficient to prove that there exists
$r>0$, such that
\begin{equation}
I(X) > 0,\forall X \in \partial {B_r}(0).
\end{equation}
We claim: $X=0$ is an isolated zero point of $I(X)$. In fact, 
for $$G(X) =\big({G_1}(X),{G_2}(X),......,{G_n}(X)\big)^T,$$
 so
$G_i(X)=G_i'(Y_i)X$, where $Y_i$ connects 0 to $X$ $(i=1,2,...,n)$.
Define a continuous function $\beta (X)$ as
\[\beta \left( X \right) = \left\{ \begin{array}{l}
{\big({G'_1}({Y_1}),{G'_2}({Y_2}),...,{G'_n}({Y_n})\big)^T}, X \ne 0,\\
~G'(0),~~~~~~~~~~~~~~~~~~~~~~~~~~~~~~~~~~~~~~~~~X = 0.
\end{array} \right.\]
Thus $G(X)=\beta (X)X$, $\forall X\in \mathbb{R}^n$.
Define
$$\gamma ({X_1},{X_2}...,{X_n}) = {\big({G'_1}({X_1}),{G'_2}({X_2}),...,{G'_n}({X_n})\big)^T}.$$
Thus $\gamma(X,X,...,X)=G'(X)$ and $\gamma(Y_1,Y_2,...,Y_n)=\beta(X)$.
Therefore, \[\det \gamma (0,0,...,0) = \det I'(0) \ne 0.\]
By the continuity of $\gamma$, there exists a positive number $r>0$, such that
\[\det \gamma ({X_1},{X_2},...,{X_n}) \ne 0,~ \mbox{for}~ (X_1,X_2,...,X_n)\in B_r(0). \]
Thus $\det \beta (X)\ne 0, ~\forall X\in \displaystyle B_{r/{\sqrt n }}(0).$
Therefore, $0$ is an isolated zero point of $I(X)$.

Let
 $\displaystyle \alpha  = \mathop {\inf }\limits_{\partial B_{r/\sqrt{n}}(0) }I(X).$
It is a positive number since $I(X)$ is continuous and nonnegative and is not zero on $\displaystyle\partial B_{r/\sqrt{n}}(0).$
Thus, $I(X)$ satisfies the condition-(a) of  Theorem 2.3.

Assume that  $c$ is a critical value of $I$, that is, $\exists~ X_c\in \mathbb{R}^n $, such that $I'(X_c)=0$. Thus
\[0 < \alpha  \le c = I({X_c}) = 0.\]
Obviously, it contradicts. By Theorem \ref{mM}, the functional $I(X)$ does not satisfy the
condition-(b) i.e., the $(PS)_c$ condition does not hold. Hence, there exists an unbounded sequence $\{X_k\} \subset \mathbb{R}^n $, such that
\[ (i)\;\; I({X_k}) \to c;~~
(ii)\;\; I'({X_k}) \to 0.~~\]
Suppose that $\{X_k\} $ is bounded in  $\mathbb{R}^n$,
then $\{X_k\} $ has a weak convergence subsequce.  It is also strong convergent in $\mathbb{R}^n$. This contradicts with the fact that $I(X)$ does not satisfy the $(PS)_c$ conditon.

\vskip0.1in

\noindent{\bf Case (1):}  $Spec(F+F^T)\subseteq (\varepsilon,+\infty).$  We let  $\mu_1$ denote the minimum eigenvalue of a Hermitian matrix $A$. That is,
\begin{equation}\label{key1}
{\mu _1} = \mathop {\inf }\limits_{Y \ne 0} \frac{{{Y^T}AY}}{{{Y^T}Y}}.
\end{equation}
Set $A = G'(X_k) +G'(X_k)^T ~\mbox{and}~ Y=G(X_k)$. By (\ref{key1}), we obtain
\begin{equation}\label{key2}
\begin{aligned}
{\mu _1}\left( {{X_k}} \right) &\le \frac{{G{{({X_k})}^T}(G'({X_k}) + G'{{({X_k})}^T})G({X_k})}}{{G{{({X_k})}^T}G({X_k})}}\\
& = \frac{{G{{({X_k})}^T}G'({X_k})G({X_k}){\rm{ + }}G{{({X_k})}^T}G'{{({X_k})}^T}G({X_k})}}{{G{{({X_k})}^T}G({X_k})}}
\\
& = \frac{{2G{{({X_k})}^T}G'({X_k})G({X_k})}}{{G{{({X_k})}^T}G({X_k})}}.
\end{aligned}
\end{equation}
By (i), one gets
\begin{equation}\label{key3}
G{({X_k})^T}G({X_k}) = I({X_k}) \to c > 0.
\end{equation}
By (ii) and (\ref{key3}), we obtain
\begin{equation}\label{key4}
\begin{aligned}
\big|2G{({X_k})^T}G'({X_k})G({X_k}) \big|&\le \left\|2 {G{{({X_k})}^T}G'({X_k})} \right\|\left\| {G({X_k})} \right\|\\
&=\left\| {I'({X_k})} \right\|{(G{({X_k})^T}G({X_k}))^{\frac{1}{2}}}\\
& = \left\| {I'({X_k})} \right\|\sqrt {I({X_k})}\to 0.
\end{aligned}
\end{equation}
Combining (\ref{key2}), (\ref{key3}) with (\ref{key4}), letting $k\to +\infty $,
we see that  $\mu_1(X_k) \leqslant  0.$  Note that  $Spec(F+F^T)\subseteq (\varepsilon,+\infty)$, we observe that 
 ${\mu _1}\left( {{X_k}} \right) \geqslant  \varepsilon.$ Letting $k \to +\infty$,  we get a contradiction.

\vskip0.12in

\noindent{\bf Case (2):}  $Spec(F+F^T)\subseteq (-\infty,-\varepsilon).$  For this case, we let $\mu _2$ be the maximum eigenvalue of a Hermitian matrix $A$. That is,
$$\displaystyle{\mu _2} = \mathop {\sup }\limits_{Y \ne 0} \frac{{{Y^T}AY}}{{{Y^T}Y}}.$$
By the same method, we obtain
$ \mu_2(X_k) \geqslant  0.$ Since $Spec(F+F^T)\subseteq (-\infty,-\varepsilon) $ and 
 $\mu _2(X_k) \leqslant  -\varepsilon$, we  also get a contradiction.
\end{proof}

\section{The proof of Theorem \ref{main2}}
\begin{proof} Firstly, we consider two cases for $Spec(F+F^T)$ in Theorem \ref{main2}.

Case (i): Suppose $Spec(F+F^T)\subseteq (-\infty,0)$ and $\exists~ M_1, M_2>0$  such that 
 for any $\lambda \in Spec(F+F^T) $, we have
 \[|\lambda |<\sum_{\lambda\in Spec(F+F^T) }|\lambda |=-\sum_{\lambda\in Spec(F+F^T) }\lambda <M_1.\]
Therefore,
\[ M _1^{n - 1}|\lambda|\geqslant \prod_{\lambda\in Spec(F+F^T)} |\lambda| >M_2>0.\]
Hence
$ \displaystyle |\lambda| >  \frac{M_2}{M_1^{n-1}}, $
i.e., $\displaystyle \lambda< -\frac{M_2}{M_1^{n-1}}. $
It implies that $Spec(F+F^T)$ can not tend to zero. 
This also contradicts to $\lambda \geqslant \varepsilon $ by the case (1) in Theorem \ref{main2}.
Thus the conclusion  holds.

\vskip0.1in

Case (ii): Suppose $Spec(F+F^T)\subseteq (0,+\infty)$ and $\exists~ M_1, M_2>0$ such that 
for any $\lambda \in Spec(F+F^T) $, we have
 \[\lambda<\sum_{\lambda\in Spec(F+F^T)}\lambda<M_1.\]
Therefore,
\[ M _1^{n - 1}\lambda\geqslant \prod_{\lambda\in Spec(F+F^T)} \lambda >M_2>0.\]
Hence, $ \displaystyle \lambda >  \frac{M_2}{M_1^{n-1}}.$
It also contradicts to $\lambda\leqslant  -\varepsilon$ in the case (2) of Theorem \ref{main2}.
 
\end{proof}

\end{CJK*}
 \end{document}